\documentclass[12pt]{amsart}
\usepackage{amsmath, amssymb, latexsym, amsthm}
\usepackage[dvips]{graphicx}

\long\def\forget#1\forgotten{}

\newtheorem{theorem}{Theorem}[section]
\newtheorem{lemma}[theorem]{Lemma}
\newtheorem{corollary}[theorem]{Corollary}
\theoremstyle{definition}
\newtheorem{definition}[theorem]{Definition}
\newcommand{\x}{\times}
\newcommand{\Q}{\mathbb{Q}}
\newcommand{\R}{\mathbb{R}}
\newcommand{\Z}{\mathbb{Z}}
\newcommand{\iso}{\cong}
\newcommand{\inv}{^{-1}}
\newcommand{\op}{\operatorname}
\newcommand{\Aut}{\op{Aut}}
\newcommand{\Supp}{\op{Supp}}
\newcommand{\bbF}{\mathbb{F}}
\newcommand{\nin}{{\not\in}}
\newcommand{\Union}{\bigcup}
\newcommand{\be}{\begin{enumerate}}
\newcommand{\ee}{\end{enumerate}}
\newcommand{\bi}{\begin{itemize}}
\newcommand{\ei}{\end{itemize}}
\renewcommand{\i}{\item}

\title[The conjugacy problem in $\ell$-groups]{The conjugacy problem and related problems in
lattice-ordered groups}
\author{W.\ Charles Holland}
\address{Department of Mathematics, Bowling Green
State University, Bowling Green, Ohio, USA}
\email{chollan@bgnet.bgsu.edu}
\author{Boaz Tsaban}
\address{Department of Applied Mathematics and Computer Science,
Weizmann Institute of Science, Rehovot 76100, Israel}
\email{boaz.tsaban@weizmann.ac.il}
\urladdr{http://www.cs.biu.ac.il/\~{}tsaban}



\begin{document}

\begin{abstract}
We study, from a constructive computational point of
view,
the techniques used to solve the conjugacy problem
in the ``generic'' lattice-ordered group $\Aut(\R)$.
We use these techniques in order
to show that for all $f,g\in\Aut(\R)$,
the equation $xfx=g$ is effectively solvable in
$\Aut(\R)$.
\end{abstract}

\maketitle

\section{Introduction}

\subsection*{The conjugacy problem}
Elements $g_1$ and $g_2$ in a group $G$ are
\emph{conjugate}
if there exists $h\in G$ such that $g_1=h\inv g_2 h$.
The \emph{conjugacy problem} for a given group $G$ is
the question whether there exists an effective
procedure
to determine whether $g_1$ and $g_2$ are conjugate,
given arbitrary $g_1,g_2\in G$.

This problem is of significant theoretical interest,
but recently
it became extremely important from a practical point
of view.
In \cite{Anshel} and \cite{KoLee}, a family of
cryptosystems
was suggested, whose strength heavily depends on the
intractability of variants of the conjugacy problem
in the underlying group.
It seemed that for $G=B_n$, the \emph{Braid group}
with $n$ strands,
the goal of achieving a secure cryptosystem was
reached,
but recent results \cite{GKTTV06, GKTTV05, Dennis} suggest that
$B_n$ is not a good candidate and the search for a
better group
has revived.

\subsection*{Lattice-ordered groups.}
A \emph{partially ordered group} is a group $G$
endowed with a
partial ordering $\le$ which is respected by the group
operations, that
is, for each $g_1,g_2\in G$ such that $g_1\le g_2$,
$xg_1\le xg_2$
and $g_1x\le g_2x$ for all $x\in G$.
If the underlying partial order $\le$ on $G$ is a
lattice
(that is, for each $g_1,g_2\in G$ there exists a least
upper bound
$g_1\lor g_2\in G$ and a greatest lower bound
$g_1\land g_2\in G$),
then we say that $G$ is a \emph{lattice-ordered
group}.

If $G$ is a lattice-ordered group, then
the lattice operations distribute over each other, and
the group operation distributes over the
lattice operations, too.
Consequently, any element in a lattice-ordered group
generated by $\{x_1,\dots,x_n\}$ can be
written in the form
$$w(x_1,\ldots,x_n)=\bigwedge_i\bigvee_ju_{ij}(x_1,\ldots,x_n)$$
where each expression $u_{ij}(x_1,\ldots,x_n)$ is an
element of the free group on
$\{x_1,\ldots,x_n\}$.
The form above for $w(x_1,\ldots,x_n)$ is not unique.
In \cite{HM}, an algorithm was given to determine
whether two given expressions of this form
represent the same element of the \emph{free
lattice-ordered group}
$\bbF_n$ (and, therefore, the same element in
\emph{every}
lattice-ordered group $G$).

Let $\Aut(\R)$ denote the collection
of all order preserving bijections $f:\R\to\R$, that
is, order automorphisms of $\R$.
Observe that each $f\in\Aut(\R)$ is continuous.
$\Aut(\R)$, with the operation
of composition, is a group which is lattice-ordered
by:
$$f\le g\quad\mbox{if}\qquad f(x)\le g(x)\mbox{ for
all }x\in \R.$$
The lattice operations are defined by
\begin{eqnarray*}
(f\lor  g)(x) & = & \max\{f(x),g(x)\}\\
(f\land g)(x) & = & \min\{f(x),g(x)\}
\end{eqnarray*}
for each $x\in\R$.
In 1963, Holland proved that
every lattice-ordered group can be
embedded in the lattice-ordered group
$\Aut(\Omega,\le)$ of automorphisms of a totally
ordered set $(\Omega,\le)$. This is discussed
in detail in section 7.1 of \cite{Glass}.
A particular case of this theorem is, that
the free lattice-ordered group $\bbF_n$
can be embedded in $\Aut(\R)$.
Consequently, $\Aut(\R)$ satisfies a given equation
$w(x_1,\ldots,x_n)=u(x_1,\ldots,x_n)$ if, and only if,
every lattice-ordered group satisfies
this equation.

\subsection*{Parametric equations}
Because of the generic nature of the lattice-ordered
group $\Aut(\R)$,
it would be interesting to know which elements of this
group are conjugate.
A simple conjugacy criterion was given in \cite{H}.
In Section 2 we analyze this treatment of the
\emph{conjugacy problem in $\Aut(\R)$} from a
computational point of view,
and show that in fact, there exists an effective
definition of the
conjugator when the given elements are conjugate.
The conjugacy problem in $\Aut(\R)$
is a specific case of an \emph{equation with
parameters} from $\Aut(\R)$.
Thus, a natural extension of the conjugacy problem in
this group is,
which equations with parameters in $\Aut(\R)$ have
solutions in
$\Aut(\R)$. We solve several problems of this type in
Section 3. In particular, we show that
every element of $\Aut(\R)$ is a commutator (that is,
for each
$g\in \Aut(\R)$ there exist $x,y\in \Aut(\R)$ such
that $x\inv y\inv  xy=g$),
and that the equation $xfx=g$ is effectively solvable
in $\Aut(\R)$.

\subsection*{Effectiveness}
When dealing with parametric equations in $\Aut(\R)$, we use the following natural
model of computation: The parameters appearing in the equation are
treated as ``black box'' functions, that is, the allowed operations
in our model are
evaluation of any of the parameters at any desired point in $\R$,
as well as the basic arithmetic operations (addition, subtraction,
multiplication and division), and any (finite) composition of these.

Moreover, we consider the basic arithmetic operations as computationally
negligible. Thus, in this model, a solution to a given
parametric equation (that is, well defined elements of $\Aut(\R)$ which
satisfy the equation when substituted for the variables)
is \emph{effective} if
its evaluation at each given point requires only finitely
many evaluations of the functions appearing as
parameters in the equations.

\subsubsection*{Notational convention}
For the rest of this paper, we use the convention that
the functions are evaluated from left to right,
that is, the value of $g$ at $\alpha$ is $\alpha g$
and the value of $gf$ at
$\alpha$ is $\alpha gf = (\alpha g)f$.

\section{The conjugacy problem}

For $g\in \mbox{Aut}(\R)$, let
\begin{eqnarray*}
\Supp(g)=\{\alpha\in\R\ :\  \alpha \neq \alpha g\}\\
\Supp^+(g)=\{\alpha\in\R\ :\  \alpha < \alpha g\}\\
\Supp^-(g)=\{\alpha\in\R\ :\  \alpha g < \alpha\}
\end{eqnarray*}
Then $\Supp^+(g)$ and $\Supp^-(g)$ are disjoint open
subsets of $\R$,
and $\Supp(g)=\Supp^+(g)\cup\Supp^-(g)$.
Consequently, $\Supp(g)$ is a disjoint
union of open intervals (the \emph{components} of
$\Supp(g)$),
where each interval is a component
of either $\Supp^+(g)$ (a \emph{positive component})
or of $\Supp^-(g)$
(a \emph{negative component}).

We now describe a useful method to obtain a partition
of
a component of $\Supp(g)$ into a sequence of half-open
intervals.
Suppose $I$ is a positive component of $\Supp(g)$ and
$\alpha\in I$.
Then $\alpha<\alpha g$.
As $g$ is order preserving, we have that
for all $i\in\Z$, $\alpha g^i<\alpha g^{i+1}$.
Let $I'$ be the \emph{convex hull} of $\{\alpha g^i\
:\  i\in\Z\}$, that is,
$$I' = \Union_{i\in\Z} [\alpha g^i,\alpha g^{i+1}).$$
Then $I'\subseteq I$. Moreover, $I'g=I'$.
If $I'$ has an upper bound, then it has a least upper
bound $\gamma$, and
$\lim_{n\to\infty}\alpha g^n = \gamma$. As $g$ is
continuous,
$\gamma g=\gamma$, and so $\gamma\nin I$.
A similar result holds if $I'$ has a lower bound.
Therefore,  $I'=I$.
Similarly, for each negative component $I$ of
$\Supp(g)$ and each $\alpha\in I$,
$$I = \Union_{i\in\Z} [\alpha g^{i+1},\alpha g^i).$$

The following lemma is an extension of the corresponding
lemma from \cite{H}.
Recall that if $\alpha$ lies in a positive component of
a function $g$, then $\alpha<\alpha g$, and the function
$\psi$ in the following lemma is well defined.

\begin{lemma}\label{singlebump}
Let $f,g\in\Aut(\R)$, let $I$ be a positive component
of
$\Supp(f)$ and $J$ be a positive component of
$\Supp(g)$.
Fix elements $\alpha\in I$ and $\beta\in J$.
Define the usual affine order preserving bijection
$\psi:[\alpha,\alpha g)\to[\beta,\beta f)$ by
$$\gamma\psi = \frac{\beta f-\beta}{\alpha
g-\alpha}(\gamma-\alpha)+\beta.$$
The following procedure defines an order preserving
bijection
$x:I\rightarrow J$ such that on $J$, $f=x\inv gx$,
by defining its evaluation on a given $\gamma\in I$:
\be
\i If $\gamma>\alpha$, compute $\alpha g,\alpha
g^2,\dots$ until
the first positive integer $i$ such that $\alpha
g^i\le \gamma<\alpha g^{i+1}$ is found.
\i If $\gamma<\alpha$, compute $\alpha g\inv,\alpha
g^{-2},\dots$ until
the first negative integer $i$ such that $\alpha
g^i\le \gamma<\alpha g^{i+1}$ is found.
\i Compute $\gamma x := \gamma g^{-i}\psi f^i$ by
making $i$ evaluations of $g\inv$, one evaluation
of $\psi$, and $i$ evaluations of $f$.
\ee
A similar result holds in the case that $I$ and $J$
are
negative components.
\end{lemma}

\begin{proof}
Let $\alpha\in I$, $\beta\in J$.
We may assume $\alpha<\alpha g$ and $\beta<\beta f$.
$$I = \Union_{i\in\Z} [\alpha g^i,\alpha g^{i+1})$$
and
$$J = \Union_{i\in\Z} [\beta f^i,\beta f^{i+1}).$$
Let $\psi:[\alpha,\alpha g]\rightarrow [\beta,\beta
g]$
 be the order preserving bijection defined above. For
each
 $i\in\Z$ define an order preserving bijection
$x_i:[\alpha g^i,\alpha g^{i+1})\to[\beta f^i,\beta
f^{i+1})$
by
$$x_i = g^{-i}\psi f^i,$$
and take $x = \Union_{i\in\Z}x_i$.
Then $x:I\rightarrow J$ is an
order preserving bijection, and if $\beta
f^i\le\delta<\beta f^{i+1}$,
then $\alpha g^i\le\delta x\inv<\alpha g^{i+1}$.
Therefore,
\begin{eqnarray*}
\delta x\inv gx & = & \delta x\inv gg^{-(i+1)}\psi
f^{i+1} =\\
& = & \delta x\inv g^{-i}\psi f^{i+1}=\delta x\inv
xf=\delta f.
\end{eqnarray*}
\end{proof}

The following is obvious.

\begin{lemma}\label{nobump}
 Let $I$ and $J$ be nontrivial maximal intervals
 of fixed points of $f$ and $g$, respectively.
\be
\i If $I = [\alpha_1,\alpha_2]$ and $J =
[\beta_1,\beta_2]$,
define $$\psi:[\alpha_1,\alpha_2]\rightarrow
[\beta_1,\beta_2]$$ as in Lemma \ref{singlebump};

\i If $I = (-\infty,\alpha_2]$ and $J =
(-\infty,\beta_2]$,
define $$\psi:(-\infty,\alpha_2]\rightarrow
(-\infty,\beta_2]$$ by
$\gamma\psi = \gamma - \alpha_2 + \beta_2$;

\i If $I = [\alpha_1,\infty)$ and $J =
[\beta_1,\infty)$,
define $$\psi:[\alpha_1,\infty)\rightarrow
[\beta_1,\infty)$$
by $\gamma\psi = \gamma - \alpha_1 + \beta_1$;

\i If $I = \R = J$,
define $$\psi:\R\rightarrow\R$$ by
$\gamma\psi = \gamma$.
\ee
Let $x = \psi$.
Then $x:I\rightarrow J$ is an order preserving
bijection such that on $J$,
$f = x^{-1}gx$.
\end{lemma}

The computational complexity in Lemmas
\ref{singlebump} and \ref{nobump} is
unbounded, but the procedure requires
only finitely many steps. For each given $\gamma$, the
computational complexity
of the evaluation of $\gamma x$ can be reduced from
$i$ (as defined there)
to the order of $\log_2i$ if we work in the \emph{fast
forward model},
where the computational complexity of evaluating $g^i$
and $f^i$ is independent
of $i$ (this model was studied in another context in
\cite{NaRe, FFPerms}).
In this model, step 3 of the procedure requires a
negligible amount of time,
and step 1 can be accelerated by first finding the
first $n$ such that
$\alpha g^{2^n}<\gamma<\alpha g^{2^{n+1}}$ and
continuing this binary
search in the interval $[\alpha g^{2^n},\alpha
g^{2^{n+1}})$ in a nested manner.

\begin{definition}
For an element $g\in\Aut(\R)$,
let $F(g)$ be the set of nontrivial maximal intervals of fixed
points of $g$, let $P(g)$ be the set of
positive components of $\Supp(g)$,
and let $N(g)$ be the set of negative components of
$\Supp(g)$.
The set $T(g)=P(g)\cup N(g)\cup F(g)$ inherits a total order from
$\R$. We call $T(g)$ the \emph{terrain} of $g$.
\end{definition}

Following is a simple characterization of terrains.
\begin{lemma}
Assume that $g\in\Aut(\R)$.
Give the elements of $F(g)$ the color $0$,
the elements of $P(g)$ the color $+$, and
the elements of $N(g)$ the color $-$.
Then the terrain $T(g)$ is a countable
$\{0,+,-\}$-colored totally ordered set such that
no two adjacent points are both
colored $0$.

Conversely, any countable
$\{0,+,-\}$-colored totally ordered set such that
no two adjacent points are both
colored $0$ is the terrain of some
element $g\in\Aut(\R)$.
\end{lemma}
\begin{proof}
$T(g)$ is countable because the component intervals and the
maximal nontrivial fixed point intervals are all disjoint, and each contains a
rational number. No two adjacent intervals are both fixed point
intervals, as this would contradict the maximality.

Conversely, if $T$ is a countable $\{0,+,-\}$-colored ordered set, let
$\Q$ be the set of rational numbers with the usual order and let
$S = T\x\Q$ be the lexicographically ordered product.
Then $S$ is a countable ordered set without end points, and so $S$ is
isomorphic $\Q$, and hence the Dedekind completion of $S$ is isomorphic
to the real line $\R$.
Under the isomorphism, for each $t\in T$, the Dedekind completion of
the interval $\{t\}\x\Q$ is isomorphic to an interval of the form
$\{t\}\x\R$, and we can define $g\in\Aut(\R)$ so that if
$t$ has color $+$, then $(t,x)g = (t,x+1)$, and if $t$ has color $-$,
then $(t,x)g = (t,x-1)$, and $g$ fixes all other points of the Dedekind
completion of $S$.
Then the terrain of $g$ is isomorphic to $T$.
\end{proof}

\begin{definition}
An \emph{isomorphism} of terrains $T_1$ and $T_2$ is a color- and
order-preserving bijection from $T_1$ to $T_2$. If there exists
such an isomorphism then we say that $T_1$ and $T_2$ are
isomorphic.
\end{definition}

\begin{theorem}\label{conj}
Two elements $g,f\in\Aut(\R)$ are conjugate if, and only if,
$T(f)$ is isomorphic to $T(g)$. Moreover, if an isomorphism of
$T(f)$ and $T(g)$ is given (as a ``black-box'' function), then
there exists an effective procedure defining an element
$h\in\Aut(\R)$ such that $f=h\inv gh$.
\end{theorem}
\begin{proof}
It is clear that if $C\in T(f)$ is a component of $f$,
and $h\in\Aut(\R)$, then
$Ch\in T(h^{-1}fh)$ is a component of $h^{-1}fh$ of
the same ``color''. Hence,
conjugation by $h$ induces an isomorphism of the
terrains $T(f)\cong T(h^{-1}fh)$.
Conversely, if we are given an isomorphism
$\tau:T(f)\cong T(g)$ of terrains,
then for every component $I = C\in T(f)$, and $J =
C\tau\in T(g)$, we have that
$I$ and $J$ satisfy the conditions of Lemmas
\ref{singlebump} or \ref{nobump}.
Since the union of all of the components of any
element of $\Aut(\R)$ is a dense
subset of $\R$, if $x$ is defined on each of the
intervals as in Lemmas \ref{singlebump}
 and \ref{nobump},
there is a unique extension to an element
$h\in\Aut(\R)$, and the theorem is proved.
\end{proof}

Of course, $T(f)$ may typically be infinite, but for a
large class of elements, it is finite.
For example, there are exactly three (isomorphism
classes of) one-element terrains,
and thus three conjugacy classes of the corresponding
members of $\Aut(\R)$.
There are exactly 8 two-element terrains, and so 8
conjugacy classes of corresponding
elements of Aut(R). And there are exactly 22
three-element terrains, etc.

\section{Other parametric equations}

The conjugacy problem in Section 2 can be expressed in
the following way:
Given parameters $g_1,g_2$, does there exist a $g_3\in
\Aut(\R)$ such that
$g_1^{-1}g_3^{-1}g_2g_3 = e$?
The general problem is this: given a lattice-ordered
group $G$ and an
element
$$w(x_1,\ldots,x_n)=\bigwedge_i\bigvee_ju_{ij}(x_1,\ldots,x_n)$$
of the free lattice-ordered group on
$\{x_1,\ldots,x_k,\ldots,x_n\}$,
$1\le k\le n$, and elements $g_1,\ldots,g_{k-1}\in G$,
do there exist elements
$g_k,\ldots,g_n\in G$ such that
$w(g_1,\ldots,g_{k-1},g_k,\ldots,g_n)=e$?

Another special case of this is when there is only one
parameter, and it occurs only once.
This was solved (modulo the effectiveness assertion)
in the following theorem and
corollaries in  \cite{AH}.

\begin{theorem}\label{1param}
Let $w(x_2,...,x_n)$ be a group word (not involving
the lattice operations,
and let $g\in \Aut(\R)$. Then there exists an
effective proceedure defining
$g_2,...,g_n\in \Aut(\R)$ such that $g =
w(g_2,...,g_n)$.
\end{theorem}

\begin{proof}
We define the functions $g_i$ ($i=2,\ldots,n$) on each
component
$I$ of $\Supp(g)$, and then patch the results
together.
Let $I$ be a component of $\Supp(g)$, say, a positive
component.
Choose any $\alpha\in I$.
Then, as shown in the previous section, $\alpha<\alpha
g$ and $\{\alpha g^i\}$ is
unbounded above and below in $I$.

We may write the equation
$w(x_2,\ldots,x_n)=g$ in the form
$$w(x_2,\ldots,x_n)=x_{\sigma(1)}^{\epsilon(1)}x_{\sigma(2)}^{\epsilon(2)}\cdots
x_{\sigma(m)}^{\epsilon(m)}=g,$$
where
$\sigma:\{1,\ldots,m\}\rightarrow\{2,\ldots,n\}$, and
$\epsilon(i)=\pm 1$, and we
may assume that the left-hand side is in reduced form,
that is,
$x_{\sigma(i+1)}^{\epsilon(i+1)}\not=x_{\sigma(i)}^{-\epsilon(i)}$.

Let $\cdots<\beta_i<\beta_{i+1}<\cdots$ be any
sequence of points of $I$
which has no
upper or lower bound in $I$. In each interval
$[\beta_i,\beta_{i+1})$,
choose points
$\beta_i=\gamma_{i,0}<\gamma_{i,1}<\cdots<\gamma_{i,m}=\beta_{i+1}$.

For each $\sigma(j)$ with $0<j\le m$ we can define an
order
preserving bijection $g_{\sigma(j)}\in\Aut(\R)$ such
that for each
$i\in\Z$ and each $j$,
$\gamma_{i,j-1}g_{\sigma(j)}^{\epsilon(j)}=\gamma_{i,j}$,
and
$\gamma_{\sigma(j)}$ is affine on each of the
intervals
$[\gamma_{i,k-1},\gamma_{i,k}]$.

We have that $\beta_iw(g_2,\ldots,g_n)=\beta_{i+1}$
for each $i$.
We do not necessarily have that $w(g_2,\ldots,g_n)=g$,
but we do have
that
$I$ is a positive component of $w(g_2,\ldots,g_n)$.
Therefore, by Theorem
\ref{conj}, there is (an effectively computable)
$y\in\Aut(\R)$ such that on $I$
$$w(y\inv g_2y,\ldots,y\inv g_ny)=y\inv
w(g_2,\ldots,g_n)y=g.$$
We do this on each component, letting all $x$'s be $e$
on each fixed point of $g$,
and patch the results together,
and the theorem is proved.
\end{proof}

\begin{corollary}
Every element $g\in\Aut(\R)$ is a commutator.
\end{corollary}
\begin{proof}
Take $g=x\inv y\inv xy$.
\end{proof}

\begin{corollary}
Every element $g\in\Aut(\R)$ has an $n$th root for
each positive
integer $n$.
\end{corollary}
\begin{proof}
Take $g=x^n$.
\end{proof}

Let us now consider the case of two parameters, but
only one variable.
Two special cases of this are considered in the next
theorem.

\begin{theorem}\label{1var}
Let $\epsilon(i)=\pm 1$, and consider the equation
$x^{\epsilon(1)}gx^{\epsilon(2)}f\inv=e$
in $\Aut(\R)$. Then:
\be
\i If $\epsilon(1)=-\epsilon(2)$, then the equation
has a solution if and
only if $T(f)\iso T(g)$.
\i If $\epsilon(1)=\epsilon(2)$, then the equation has
a solution for all
$f,g$.
\ee
Moreover, when these equations have solutions, they
have effectively
defined solutions.
\end{theorem}
\begin{proof}
(1) This is Theorem \ref{conj}.

(2) We write the equation in the form $xgx=f$. Since
$f\inv(fg)f=gf$, by
Theorem \ref{conj}, $T(fg)\iso T(gf)$. In particular,
if $I$ is a component of
$\Supp(fg)$,
then $If$ is the corresponding component of
$\Supp(gf)$. Suppose, for
example that
$fg$ is positive on $I$. Then $gf$ is positive on
$If$. Choose $\alpha\in
I$.
Then
$$\cdots<\alpha<\alpha fg<\alpha (fg)^2<\cdots$$
is unbounded in $I$.
It follows that
$$\cdots<\alpha f(gf)\inv=\alpha g\inv<\alpha f<\alpha
f(gf)<\cdots$$
is unbounded in $If$.
Choose $\beta \in If$ so that $\alpha
g\inv<\beta<\alpha f$. Then
$$\cdots<\alpha<\beta g<\alpha (fg)<\beta g(fg)<\alpha
(fg)^2<\beta
g(fg)^2<\cdots$$
and
$$\cdots<\beta<\alpha f<\beta (gf)<\alpha f(gf)<\beta
(gf)^2<\alpha
f(gf)^2<\cdots,$$
and each of these sequences is unbounded in the
corresponding component.

Let $\psi :[\alpha,\beta g)\rightarrow [\beta,\alpha
f)$ be any order
preserving
bijection between the real intervals, for example the
affine one. We now define an order preserving
bijection
$x:I\rightarrow If$ by extending $\psi$ in the
following way. For
$\gamma\in I$:
$$\gamma x=\left\{\begin{array}{ll}
\gamma (fg)^{-i}\psi (gf)^i, &\mbox{ if }\alpha
(fg)^i\le\gamma<\beta
g(fg)^i\\
\gamma (fg)^{-i}g\inv\psi\inv f(gf)^i, &\mbox{ if
}\beta g(fg)^i\le
\gamma<\alpha (fg)^{i+1}.
\end{array}
\right. $$
Then $x$ is, indeed, an order preserving bijection of
$I$ onto $If$. And
on $I$, $xgx=f$ because:
if $\alpha (fg)^i\le\gamma<\beta g(fg)^i$ then $\gamma
x=\gamma
(fg)^{-i}\psi (gf)^i$, and
so $$\beta (gf)^i=\alpha (fg)^i(fg)^{-i}\psi
(gf)^i\le\gamma x<\beta
g(fg)^i(fg)^{-i}\psi (gf)^i
=\alpha f(gf)^i$$
and so
$$\beta g(fg)^i=\beta (gf)^ig\le \gamma xg<\alpha
f(gf)^ig=\alpha
(fg)^{i+1}$$
which implies
\begin{eqnarray*}
\gamma xgx&=&(\gamma x)g(fg)^{-i}g\inv\psi\inv
f(gf)^i\\
&=&
(\gamma (fg)^{-i}\psi (gf)^i)g(fg)^{-i}g\inv\psi\inv
f(gf)^i\\
&=&
\gamma (fg)^{-i}(fg)^if\\
&=&
\gamma f;
\end{eqnarray*}
and in the other case,
$\beta g(fg)^i\le\gamma<\alpha (fg)^{i+1}$, so $\gamma
x=\gamma
(fg)^{-i}g\inv\psi\inv f(gf)^i$,
whence
\begin{eqnarray*}
\alpha f(gf)^i&=&\beta g(fg)^i(fg)^{-i}g\inv\psi\inv
f(gf)^i\\
  &\le&\gamma x\\
  &<&\alpha (fg)^{i+1}(fg)^{-i}g\inv
\psi\inv f(gf)^i\\
&=&\beta (gf)^{i+1},
\end{eqnarray*}
from which follows
$\alpha (fg)^{i+1}\le \gamma xg<\beta g(fg)^{i+1}$,
and hence
\begin{eqnarray*}
\gamma xgx&=&(\gamma x)g(fg)^{-(i+1)}\psi (gf)^{i+1}\\
&=&(\gamma (fg)^{-i}g\inv\psi\inv
f(gf)^i)g(fg)^{-(i+1)}\psi (gf)^{i+1}\\
&=&\gamma (fg)^{-i}g\inv(gf)^{i+1}\\
&=&\gamma f.
\end{eqnarray*}

Repeating this process on each of the components of
$\Supp(fg)$
(and defining $x=f$ on the fixed points of $fg$),
produces an
$x\in\Aut(\R)$
such that $xgx=f$.
\end{proof}


\begin{thebibliography}{99}

\bibitem{AH}
Samson Adeleke and W.\ C.\ Holland,
\emph{Representation of order automorphisms by words},
Forum Mathematicum \textbf{6} (1994),
315--321.

\bibitem{Anshel}
I.\ Anshel, M.\ Anshel, and D.\ Goldfeld,
\emph{An algebraic method for Public-Key
Cryptography},
Mathematical Research Letters \textbf{6} (1999),
287--291.

\bibitem{GKTTV05}
D.\ Garber, S.\ Kaplan, M.\ Teicher, B.\ Tsaban, and U.\ Vishne,
\emph{Probabilistic solutions of equations in the braid group},
Advances in Applied Mathematics \textbf{35} (2005),
323--334.

\bibitem{GKTTV06}
\bibitem{LBCS}
D.\ Garber, S.\ Kaplan, M.\ Teicher, B.\ Tsaban, and U.\ Vishne,
\emph{Length-based conjugacy search in the Braid group},
Contemporary Mathematics \textbf{418} (2006), 75--87.

\bibitem{Glass}
A.\ M.\ W.\ Glass,
\emph{Partially Ordered Groups},
World Scientific, 1999.

\bibitem{Dennis}
D.\ Hofheinz and R.\ Steinwandt,
\emph{A practical attack on some Braid group based
cryptographic primitives},
International Workshop on Practice and Theory in
Public Key Cryptography, PKC 2003 Proceedings,
LNCS \textbf{2567} (2002),
187--198.

\bibitem{H}
W.\ C.\ Holland,
\emph{The lattice-ordered group of automorphisms of an
ordered set},
Michigan Math.\ J.\ \textbf{10} (1963), 399--408.

\bibitem{HM}
W.\ C.\ Holland and S.\ H.\ McCleary,
\emph{The word problem for free lattice-ordered
groups},
Houston J.\ Math.\ \textbf{5} (1979), 99--105.

\bibitem{KoLee}
K.\ H.\ Ko, S.\ J.\ Lee, J.\ H.\ Cheon, J.\ W.\ Han,
J.\ Kang, and C.\ Park,
\emph{New Public-Key Cryptosystem using Braid groups},
Advances in Cryptology -- Crypto 2000 Proceedings,
LNCS \textbf{1880},
166--183.

\bibitem{NaRe}
M.\ Naor and O.\ Reingold,
\emph{Constructing Pseudo-Random Permutations with a
Prescribed Structure},
Journal of Cryptology \textbf{15} (2002),
97--102.

\bibitem{FFPerms}
B.\ Tsaban,
\emph{Permutation graphs, fast forward permutations, and sampling the cycle structure of a permutation},
Journal of Algorithms \textbf{47} (2003),
104--121.
\end{thebibliography}
\end{document}